\documentclass[11pt]{article}

\usepackage[T1]{fontenc}
\usepackage[utf8]{inputenc}
\usepackage[english]{babel}  
\usepackage{lmodern}
\usepackage{amsmath,amssymb,amsthm,mathtools}
\usepackage{microtype}
\usepackage{enumitem}
\usepackage[a4paper,margin=30mm]{geometry}
\usepackage[hidelinks]{hyperref}

\usepackage{caption}                      
\usepackage{graphicx}  
\usepackage{float}      
\usepackage{xcolor}

\usepackage{tikz}
\usetikzlibrary{arrows.meta,calc,decorations.pathreplacing}

\usepackage{verbatim}
\usepackage[noabbrev, nameinlink, capitalise]{cleveref}
\crefdefaultlabelformat{#2\text{\textup{#1}}#3}

\makeatletter
\def\thm@space@setup{%
	\thm@preskip=10pt
	\thm@postskip=6pt
}
\makeatother

\newtheorem{theorem}{Theorem}[section]
\crefname{theorem}{Theorem}{Theorems}
\newtheorem*{theorem*}{Theorem}
\crefname{theorem*}{Theorem}{Theorems}
\newtheorem{lemma}[theorem]{Lemma}
\crefname{lemma}{Lemma}{Llemas}
\newtheorem{proposition}[theorem]{Proposition}
\newtheorem{corollary}[theorem]{Corollary}
\newtheorem{conjecture}[theorem]{Conjecture}
\crefname{conjecture}{Conjecture}{Conjectures}

\theoremstyle{definition}
\newtheorem{definition}[theorem]{Definition}
\theoremstyle{remark}
\newtheorem{remark}[theorem]{Remark}

\numberwithin{equation}{section}

\DeclareMathOperator{\esssup}{ess\,sup}
\DeclareMathOperator{\essinf}{ess\,inf}

\DeclareMathOperator{\supp}{supp}
\DeclareMathOperator{\dm}{dm}

\newcommand{\R}{\mathbb{R}}
\newcommand{\C}{\mathbb{C}}
\newcommand{\Torus}{\mathbb{T}}
\newcommand{\PW}{PW}
\newcommand{\one}{\mathbf{1}}
\newcommand{\cC}{\mathcal{C}}
\newcommand{\cA}{\mathcal{A}}
\newcommand{\cU}{\mathcal{U}}

\newcommand{\TT}{\mathcal{T}}
\newcommand{\WW}{\mathcal{W}}
\newcommand{\diamess}{\operatorname{diam}_{\mathrm{ess}}}

\usepackage{anysize}\marginsize{2.5cm}{2.5cm}{1.1cm}{1.7cm}  
\usepackage{setspace}
 
\usepackage{parskip}
\usepackage{yhmath}

\begin{document}

\thispagestyle{empty}
\vspace*{0cm}

\begin{center}
	\Large
	\textbf{Phase symmetrization and a no-aliasing \\ concentration principle in Paley--Wiener spaces} \\[12pt]
	
	\small
	\textsc{Oriol Baeza Guasch} \\[6pt]
	\textit{Universitat Politècnica de Catalunya} 
	
\end{center}

\vspace{12pt}

\begin{center}
	\begin{minipage}{0.85\textwidth}
		\begin{spacing}{1}
			\fontsize{8}{12}\selectfont
			\setlength{\parskip}{0pt}
			\setlength{\abovedisplayskip}{5pt}
			\setlength{\belowdisplayskip}{5pt}
			\setlength{\abovedisplayshortskip}{8pt}
			\setlength{\belowdisplayshortskip}{8pt}
			
			\textsc{Abstract.} We study a constructive version of the Donoho--Stark concentration problem in a no-aliasing regime. Let \(f\in\PW(\WW)\), where \(\WW\subset\mathbb R\) is an interval of length \(W\), and let \(\TT\subset\mathbb R\) have measure \(T\) and essential diameter \(D\). If \(WD\leq1\), then the phase-aligned function \(g\), defined by \(\widehat g=|\widehat f|\), satisfies
			\[
			\int_\TT |f(t)|^2\,dt
			\leq
			\int_{-T/2}^{T/2}|g(t)|^2\,dt.
			\]
			Thus each band-limited function is associated with an explicit comparator having the same Fourier modulus and \(L^2\)-norm. The proof combines two independent principles: Fourier-phase alignment maximizes concentration on a centered interval for fixed Fourier modulus, while a one-period geometric estimate controls the Fourier transform of the observation set. We show that the latter mechanism is essentially sharp: beyond the one-period threshold, separated components may alias under the exponential map and reinforce coherently. Analogous phase and coefficient rearrangement results are obtained for polynomials on the unit circle, including a constant-one Montgomery-type inequality in the corresponding no-aliasing regime. The results identify both the scope and the limitation of this modewise approach to concentration.
	\end{spacing}
	\end{minipage}
\end{center}

\vspace{14mm}

\tableofcontents

\newpage

\section{Introduction}
\label{sec:introduction}

\subsection{The operator-level concentration conjecture} 

We begin by recalling the original conjecture by Donoho and Stark. For that purpose, we introduce the following operators, defined for any $f\in L^2(\mathbb{R})$. 

First, the time-limiting operator for a given measurable subset $\TT$ is
\begin{equation}
	(P_\TT f)(t) = \begin{cases} f(t) & t\in \TT, \\ 0 &\text{otherwise}, \end{cases} 
\end{equation}
and second, the frequency-limiting operator for a given measurable subset $\WW$ is
\begin{equation}
	(P_\WW f)(t) = \int_\WW e^{2\pi i w t} \widehat{f}(w) \, dw,
\end{equation}
where $\widehat f$ is the Fourier transform of $f$, with the convention
\begin{equation}
	\widehat f (w) =\int_\mathbb{R} f(t) e^{-2\pi i w t} \, dt.
\end{equation}

Note that the Paley--Wiener space for $\WW\subset \R$ a bounded interval is precisely
\begin{equation}
	\PW(\WW) =
	\bigl\{
	f\in L^2(\R):
	\supp\widehat f\subseteq\WW
	\bigr\} = 
	\bigl\{
	f\in L^2(\R):
	f = P_\WW f
	\bigr\}
	.
\end{equation}

Then, their conjecture can be stated as follows.
\begin{conjecture}[Donoho and Stark, Conj. 1~\cite{DonohoStark1989}] \label{conj:continuous}
	The supremum $\sup \|P_\WW P_\TT\|$, where $\WW$ is a fixed bounded interval and $\TT$ ranges over measurable subsets with fixed measure, is attained when $\TT$ is also an interval. 
\end{conjecture}

However, we shall use instead the equivalent adjoint formulation, where we consider $\sup \|P_\TT P_\WW \|$ instead. This formulation is more convenient for the concentration interpretation used below, while being equivalent to the original conjecture. Indeed, note that both operators $P_\TT$ and $P_\WW$ are self-adjoint orthogonal projections, so $(P_\TT P_\WW)^* = P_\WW P_\TT$, hence the norms of the operators are equal $\|P_\TT P_\WW\|=\|P_\WW P_\TT\|$.

To introduce the interpretation of this adjoint formulation, we begin by defining the concentration of $f\in L^2(\mathbb{R})$ in a measurable set $\TT \subseteq \mathbb{R}$ by
\begin{equation}
	c_\TT (f) = \dfrac{\displaystyle\int_{\TT} |f(t)|^2 \, dt}{\displaystyle\int_{\mathbb{R}} |f(t)|^2 \, dt}
\end{equation}
and the optimal band-limited concentration on $\TT$ by
\begin{equation}
	\mathcal{C}_\WW(\TT) = \sup_{\substack{0\neq f\in L^2 \\ f = P_\WW f}} c_\TT (f)
\end{equation}
Now, because $\|P_\WW\| =1$ and $P_\WW$ is the orthogonal projection of $L^2(\R)$ onto $\PW (\WW)$ it follows straightforwardly that

\begin{equation}
	\|P_\TT P_\WW\|^2 = \sup_{\substack{0\neq f\in L^2 \\ f = P_\WW f}} \frac{\|P_\TT f\| ^2}{\| f\| ^2} = \sup_{\substack{0\neq f\in L^2 \\ f = P_\WW f}} \dfrac{\displaystyle\int_\TT |f(t)|^2 \, dt }{\displaystyle\int_\mathbb{R} |f(t)|^2 \, dt } = \sup_{\substack{0\neq f\in L^2 \\ f = P_\WW f}} c_\TT(f) = \mathcal{C}_\WW(\TT),
\end{equation}

which we might interpret as calculating the concentration in the measurable set $\TT$, and restricting our attention to functions whose Fourier transform are supported in $\WW$. 

Altogether, in terms of the optimal concentration constant, Conjecture~\ref{conj:continuous} asserts that for $\WW$ an interval, $\TT$ any measurable subset and $I$ an interval of size $|I| = |\TT|$, it holds that
\begin{equation}
	\mathcal{C}_\WW (\TT) \leq \mathcal{C}_\WW (I).
\end{equation}

A very important remark here is that the inequality is an operator-level extremal statement: the maximizing function on \(I\) is allowed to differ completely from a given function on \(\TT\).

Finally, it is also worth mentioning that when both sets $\WW$ and $\TT$ are intervals, the concentration operator has been widely studied. In particular, the eigenfunctions of the associated concentration operator are the Prolate Spheroidal Wave Functions, with their characterization described in detail in the work of Slepian, Pollak and Landau \cite{SlepianPollak1961, Slepian1978, LandauPollak1961}, among others.


\subsection{Constructive function-by-function comparisons}

In \cite{DonohoStark1993}, Donoho and Stark proved the conjecture in the regime $|\WW||\TT|\leq 0.8$ by working in its original time-limited formulation. Given a function $f$ supported on $\TT$, they replace $|f|$ by its symmetric decreasing rearrangement $|f|^*$, which is supported on the centered interval of measure $|\TT|$, and prove that the corresponding concentration in the frequency interval $\WW$ increases.

\begin{definition}
	Given a function $f$, its symmetric decreasing rearrangement $|f|^*$ is the unique nonnegative, even function that is non-increasing as $|t|$ increases and has the same measure of its level sets,
	\begin{equation}
		\left\{ t \colon |f|^*(t) > \alpha\right\} = \left(- \frac{\mu_f(\alpha)}{2}, \frac{\mu_f(\alpha)}{2} \right), \qquad \mu_f(\alpha) = |\{|f| > \alpha \}|.
	\end{equation}
	In particular, if $f$ has essential support of finite measure $T$, then $|f|^*$ is essentially supported on $[-T/2,T/2]$.
\end{definition}

This transformation allows them to use the Riesz rearrangement inequality, in the following form.
\begin{lemma}[Hardy, Littlwood and P\'olya, Thm. 380~\cite{hardy1934inequalities}] \label{lm:hlp}
	Let $f,g$ and $h$ be positive functions. Then,
	\[
	\int_\mathbb{R} \int_\mathbb{R} f(x) g(y) h(x-y) \, dx \, dy \, \leq \int_\mathbb{R} \int_\mathbb{R} f^*(x) g^*(y) h^*(x-y) \, dx \, dy.
	\]
\end{lemma}

Thus, their argument is function-by-function in the original formulation. The equivalence of the two operator norms transfers the resulting optimal inequality to the adjoint formulation, but it does not transfer the explicit rearranged comparator to each individual band-limited function.

In a similar spirit, the present work provides a new function-by-function comparison, now in the adjoint formulation and through a different explicit transformation. That is, instead of rearranging $|f|$ in the physical variable in the original time-limited formulation, we remove the phase of $\widehat f$ in the adjoint band-limited formulation. 

Both approaches therefore establish a prescribed function-by-function comparison in their respective formulations, which is stronger than the corresponding operator-level extremal statement, but each requires an additional restriction, governed by different mechanisms.

In their case, the numerical restriction $|\WW||\TT| \leq 0.8$ enters after applying Lemma~\ref{lm:hlp} because, on the range of differences relevant to the quadratic form, the sinc kernel must agree with its symmetric decreasing rearrangement.

In our case, the numerical restriction will be 
\begin{equation}
	|\WW| \diamess (\TT) \leq 1,
\end{equation}

where recall that the essential diameter is defined as
\begin{equation}\label{eq:essential-diameter}
	\diamess(\TT)
	=
	\esssup \TT-\essinf \TT
	\in[0,\infty].	
\end{equation}

However, our condition is not solely a time-bandwidth restriction: it also controls the spatial separation of the components of $\TT$, and arises from the requirement that their exponential images do not alias.


\subsection{Main result and scope}

Our principal result is the following.

\begin{theorem}[Constructive no-aliasing concentration]
\label{thm:main-introduction}
Let \(\WW\subset\R\) be an interval of positive length \(W = |\WW|\), and let \(\TT \subset\R\) be measurable with positive finite measure \(T = |\TT|\) and finite diameter $D=\diamess(\TT)$. Assume 
\begin{equation}
    W D\leq1.
    \label{eq:main-condition-intro}
\end{equation}
For \(f\in\PW(\WW)\), define \(g\in\PW(\WW)\) by
\begin{equation}
	\widehat g=|\widehat f|.
	\label{eq:g-main-section}
\end{equation}
Then, for \(I_T = [-T/2, T/2]\), it holds
\begin{equation}
	\int_\TT|f(t)|^2\,dt
	\leq
	\int_{I_T}|g(t)|^2\,dt.
	\label{eq:function-by-function-main}
\end{equation}
Moreover,
\(
\|g\|_2=\|f\|_2
\),
and hence
\begin{equation}
	c_\TT(f)
	\leq
	c_{I_T}(g).
	\label{eq:normalized-function-by-function}
\end{equation}
\end{theorem}

The constructive character of this result should be distinguished from unrestricted interval optimization.  We do not claim that \(g\) is the global maximizer of concentration on \(I_T\), but rather the optimizer among all $h\in \PW(\WW)$ satisfying $|\widehat h| = |\widehat f|$ almost everywhere.

Also, since \(T \leq D\), condition \eqref{eq:main-condition-intro} implies \(W T \leq1\), while the converse is false: a set may have very small measure and arbitrarily large diameter. Consequently, the geometric regime considered here and the time--bandwidth regime of Donoho and Stark are not nested: our hypothesis permits, for instance, intervals with $0.8<W T \leq1$, whereas their condition $W  T \leq0.8$ permits sets of arbitrarily large diameter. Hence neither theorem subsumes the other. Moreover, on the overlap of the two regimes, they provide different prescribed comparators in different formulations.

The proof of Theorem~\ref{thm:main-introduction} separates into two independent principles.

\begin{enumerate}[label=\textup{(\roman*)}]
\item
For the centered interval \(I_T\), removal of the Fourier phase maximizes concentration among all functions having the same Fourier modulus, provided \(W T\leq1\).

\item
For a measurable set \(\TT\), the pointwise estimate
\[
    |\widehat{\one_\TT}(\xi)|
    \leq
    \widehat{\one_{I_T}}(\xi) = \frac{\sin\left(\pi |\TT| \xi\right)}{\pi \xi}
\]
holds whenever
\(
  |\xi|\diamess(\TT)\leq1
\), where the value at $\xi = 0$ is understood by continuity and equals $T$.
\end{enumerate}

For \(\omega,\eta\in\WW\), one has \(|\eta-\omega|\leq W\). Hence the condition \(W	 D\leq1\) places every difference frequency occurring in the quadratic form within the no-aliasing regime, which allows the two principles to be combined.

The diameter condition is not merely a technical assumption. The map
\[
    t\longmapsto e^{2\pi i\xi t}
\]
is periodic with period \(1/|\xi|\).  If a set spans more than one period, components which are disjoint on \(\R\) may overlap after projection onto the circle. Their Fourier contributions can then align rather than cancel. Indeed, we provide a two-component example showing that the pointwise interval estimate may fail when
\(
 |\xi|\diamess(\TT)>1
\), even by an arbitrarily small amount.

The failure of this pointwise comparison beyond the no-aliasing regime does not imply failure of the operator-level conjecture; this distinction and the limitations of the present method are discussed further in Section~\ref{sec:limitations}.

\paragraph{Note on the revised version.}
An earlier version of this manuscript used a geometric lemma asserting a pointwise estimate depending only on the total measure of a finite union of intervals.  That lemma is false because the exponential projection may identify separated components modulo one period.  The present version replaces it by the no-aliasing hypothesis
\(
 |\xi|\diamess(\TT)\leq1
\), corrects all dependent statements, and reorganizes the paper around the surviving constructive concentration principle. The aliasing counterexample and the limitations of the modewise method are recorded explicitly in Sections~\ref{sec:aliasing} and
\ref{sec:limitations}.

The organization of the paper is as follows. Section \ref{sec:phase-optimization} establishes fixed-modulus phase optimization, while Sections~\ref{sec:no-aliasing} and \ref{sec:aliasing} prove the no-aliasing estimate and analyze its sharpness. Section \ref{sec:constructive-concentration} combines these ingredients to prove the main Paley–Wiener theorem. Sections~\ref{sec:polynomial} and \ref{sec:coefficient-rearrangement} develop the polynomial and coefficient-rearrangement analogues. The final section discusses the limitations of the modewise method and its relation to the unrestricted conjecture.

\section{Phase optimization on a centered interval}
\label{sec:phase-optimization}

This section isolates the part of the argument that depends only on the
spectral phase.  Let \(\WW\subset\R\) be an interval of length \(W = |\WW|\),
and let \[
    I_T=[-T/2,T/2].
\]
The Fourier transform of its indicator is
\begin{equation}
    \widehat{\one_{I_T}}(\xi)
    =
    \begin{cases}
    \displaystyle
    \frac{\sin(\pi T\xi)}{\pi\xi}&\xi\neq0,\\[2mm]
    T&\xi=0.
    \end{cases}
    \label{eq:interval-indicator-fourier}
\end{equation}
If \(WT\leq1\), then
\begin{equation}
    \widehat{\one_{I_T}}(\xi)\geq0
    \qquad
    (|\xi|\leq W).
    \label{eq:interval-kernel-positive}
\end{equation}

\begin{lemma}[Quadratic-form representation]
\label{lem:quadratic-form-interval}
Let \(f\in\PW(\WW)\).  Then
\begin{equation}
    \int_{I_T}|f(t)|^2\,dt
    =
    \iint_{\WW\times\WW}
       \widehat f(\omega)
       \overline{\widehat f(\eta)}
       \widehat{\one_{I_T}}(\eta-\omega)
       \,d\omega\,d\eta.
    \label{eq:quadratic-form-interval}
\end{equation}
\end{lemma}

\begin{proof}
Because \(\WW\) has finite measure,
\(
  \widehat f\in L^2(\WW)\subset L^1(\WW)
\).
Fourier inversion and Fubini's theorem give
\[
\begin{aligned}
    \int_{I_T}|f(t)|^2\,dt
    &=
    \int_{I_T}
       \int_\WW\widehat f(\omega)e^{2\pi i\omega t}\,d\omega
       \int_\WW\overline{\widehat f(\eta)}
                    e^{-2\pi i\eta t}\,d\eta\,dt\\
    &=
    \iint_{\WW\times\WW}
       \widehat f(\omega)\overline{\widehat f(\eta)}
       \left(
          \int_{I_T}e^{-2\pi i(\eta-\omega)t}\,dt
       \right)
       d\omega\,d\eta,
\end{aligned}
\]
which is \eqref{eq:quadratic-form-interval}.
\end{proof}

For a nonnegative function \(h\in L^2(\WW)\), define the
fixed-modulus class
\begin{equation}
    \cA(h)
    =
    \bigl\{
       g\in\PW(\WW):
       |\widehat g|=h\ \text{a.e.}
    \bigr\},
    \label{eq:fixed-modulus-class}
\end{equation}
which is the class of functions for which their Fourier transform has prescribed amplitude almost everywhere.

The following result is the Fourier-phase counterpart of the symmetric decreasing rearrangement used by Donoho and Stark: both identify an explicit representative of a prescribed equimeasurable or fixed-modulus class that improves concentration.

\begin{theorem}[Phase alignment for fixed Fourier modulus]
\label{thm:fixed-modulus-phase}
Assume \(WT\leq1\).  Let \(h\in L^2(\WW)\), \(h\geq0\), and define
\(g\in\PW(\WW)\) by
\[
    \widehat g=h.
\]
Then
\begin{equation}
    \int_{I_T}|f(t)|^2\,dt
    \leq
    \int_{I_T}|g(t)|^2\,dt
    \qquad
    \text{for every }f\in\cA(h).
    \label{eq:fixed-modulus-phase}
\end{equation}
Thus constant spectral phase is an explicit maximizer of concentration
on the centered interval within the class \(\cA(h)\).
\end{theorem}

\begin{proof}
Write
\[
    \widehat f(\omega)=A(\omega)e^{i\phi(\omega)}
\]
on the set where \(A\neq0\).  By Lemma~\ref{lem:quadratic-form-interval},
\[
\begin{aligned}
    \int_{I_T}|f(t)|^2\,dt
    &=
    \iint_{\WW\times\WW}
       A(\omega)A(\eta)
       e^{i(\phi(\omega)-\phi(\eta))}
       \widehat{\one_{I_T}}(\eta-\omega)
       \,d\omega\,d\eta\\
    &=
    \iint_{\WW\times\WW}
       A(\omega)A(\eta)
       \cos\bigl(\phi(\omega)-\phi(\eta)\bigr)
       \widehat{\one_{I_T}}(\eta-\omega)
       \,d\omega\,d\eta.
\end{aligned}
\]
The imaginary part vanishes by symmetry under interchange of \(\omega\) and \(\eta\).  Since
\(
 \cos(\phi(\omega)-\phi(\eta))\leq1
\)
and the other factors are nonnegative by \eqref{eq:interval-kernel-positive},
\[
    \int_{I_T}|f(t)|^2\,dt
    \leq
    \iint_{\WW\times\WW}
       A(\omega)A(\eta)
       \widehat{\one_{I_T}}(\eta-\omega)
       \,d\omega\,d\eta
    =  \int_{I_T}|g(t)|^2\,dt,
\]
where the last equality follows by the definition of $g$ and Lemma~\ref{lem:quadratic-form-interval}.
\end{proof}

\begin{corollary}[Canonical phase symmetrization]
\label{cor:canonical-phase-symmetrization}
Let \(f\in\PW(\WW)\), assume \(WT\leq1\), and define \(g\) by
\(
 \widehat g=|\widehat f|
\).
Then
\[
    \|g\|_2=\|f\|_2
\]
and
\[
    c_{I_T}(f)\leq c_{I_T}(g).
\]
\end{corollary}

\begin{proof}
The norm identity follows from Plancherel's theorem.  Apply
Theorem~\ref{thm:fixed-modulus-phase} with
\(
 h=|\widehat f|
\).
\end{proof}

\begin{remark}[Translated intervals]
For the interval
\(
 I_T+t_0
\),
the fixed-modulus optimizer has Fourier transform
\[
    A(\omega)e^{-2\pi i t_0\omega}.
\]
The centered interval is distinguished only by the fact that its
canonical optimizer has constant Fourier phase.
\end{remark}

\section{No-aliasing geometry on the circle}
\label{sec:no-aliasing}

The geometric ingredient is an extremal property of the first Fourier
coefficient of a subset of one circle turn.

\begin{lemma}[Extremality of a single arc]
\label{lem:single-arc}
Let \(\cU\subset\R\) be measurable and essentially contained in an
interval of length at most \(2\pi\).  Put \(U=|\cU|\).  Then
\begin{equation}
    \left|\int_\cU e^{is}\,ds\right|
    \leq
    2\sin\left(\frac U 2\right).
    \label{eq:single-arc}
\end{equation}\
Moreover, equality is attained when \(\cU\) is an interval of length $U$, up to translation and null sets.
\end{lemma}

\begin{proof}
Since \(0\leq U\leq2\pi\),
\[
\begin{aligned}
    \left|\int_\cU e^{is}\,ds\right|
    &=
    \sup_{\theta\in\R}
    \Re\left(
      e^{-i\theta}\int_\cU e^{is}\,ds
    \right)\\
    &=
    \sup_{\theta\in\R}
    \int_\cU \cos(s-\theta)\,ds.
\end{aligned}
\]
For fixed \(\theta\), the bathtub principle says that the integral of
\(\cos(s-\theta)\) over a set of measure \(U\) is maximized by taking
the points where this function is largest.  Modulo \(2\pi\), these
points form the arc
\(
 [\theta-U/2,\theta+U/2]
\).
Therefore
\[
    \int_\cU\cos(s-\theta)\,ds
    \leq
    \int_{-U/2}^{U/2}\cos u\,du
    =
    2\sin\left(\frac U 2 \right).
\]
Taking the supremum over \(\theta\) proves the result.
\end{proof}

\begin{remark}[Chord interpretation]
If \(\cU\) is a finite union of disjoint intervals, then
\(
 i\int_\cU e^{is}\,ds
\)
is the sum of the oriented chord vectors associated with those arcs.
Lemma~\ref{lem:single-arc} says that, provided all arcs are represented
inside one turn of the circle, their resultant is no larger than the
chord generated by one arc of the same total angular length.
\end{remark}

We now transfer the circular statement to the Fourier transform of an
indicator.

\begin{theorem}[Pointwise no-aliasing estimate]
\label{thm:indicator-no-aliasing}
Let \(\TT\subset\R\) be measurable with positive finite measure \(T = |\TT|\) and finite diameter \(D=\diamess(\TT) \). If
\begin{equation}
    |\xi|D\leq1,
    \label{eq:frequency-diameter}
\end{equation}
then
\begin{equation}
    |\widehat{\one_\TT}(\xi)|
    \leq
    \begin{cases}
    \displaystyle
    \frac{\sin(\pi|\xi|T)}{\pi|\xi|}&\xi\neq0,\\[2mm]
    T&\xi=0.
    \end{cases}
    \label{eq:indicator-no-aliasing}
\end{equation}
Equivalently, for \(I_T=[-T/2,T/2]\),
\begin{equation}
    |\widehat{\one_\TT}(\xi)|
    \leq
    \widehat{\one_{I_T}}(\xi).
    \label{eq:indicator-interval-domination}
\end{equation}
\end{theorem}

\begin{proof}
The assertion at \(\xi=0\) is immediate.  Assume \(\xi\neq0\).  Up to a
null set,
\[
    \TT\subset[a,a+D],
    \qquad
    a=\essinf \TT.
\]
Set
\[
    \cU=2\pi|\xi|(\TT-a).
\]
Then \(\cU\) is essentially contained in an interval of length \(    2\pi|\xi|D\leq2\pi\), and its measure is \(|\cU| =2\pi|\xi|T\). The change of variables
\(
 s=2\pi|\xi|(t-a)
\)
gives
\[
\begin{aligned}
    |\widehat{\one_\TT}(\xi)|
    &=
    \left|
      \int_\TT e^{-2\pi i\xi t}\,dt
    \right|\\
    &=
    \frac{1}{2\pi|\xi|}
    \left|
      \int_\cU e^{-i\operatorname{sgn}(\xi)s}\,ds
    \right|.
\end{aligned}
\]
By Lemma~\ref{lem:single-arc},
\[
    |\widehat{\one_\TT}(\xi)|
    \leq
    \frac{2\sin(\pi|\xi|T)}{2\pi|\xi|}
    =
    \frac{\sin(\pi|\xi|T)}{\pi|\xi|}.
\]
Since \(T\leq D\), condition \eqref{eq:frequency-diameter} also implies
\(
 |\xi|T\leq1
\),
so the right-hand side is nonnegative.  By
\eqref{eq:interval-indicator-fourier}, it is precisely
\(
 \widehat{\one_{I_T}}(\xi)
\).
\end{proof}

For finite unions of intervals, the estimate can be written directly as
a chord-vector inequality.

\begin{corollary}[Chord-vector inequality in one period]
\label{cor:chord-vector}
Let
\[
    \TT=\bigcup_{p=1}^{r}(a_p,b_p)
\]
be a finite disjoint union ordered so that \(a_1\leq b_1\leq a_2\leq\cdots\leq a_r\leq b_r\). Put \(T=\sum_{p=1}^r(b_p-a_p)\) for its measure, and \(D =b_r - a_1\) for its diameter. If
\[
    |\xi|D\leq1,
\]
then
\begin{equation}
    \left|
      \sum_{p=1}^r
      \left(
        e^{-2\pi i\xi b_p}
        -e^{-2\pi i\xi a_p}
      \right)
    \right|
    \leq
    2\sin(\pi|\xi|T).
    \label{eq:chord-vector}
\end{equation}
\end{corollary}

\begin{proof}
For \(\xi\neq0\), direct integration yields
\[
    2\pi i\xi\widehat{\one_\TT}(\xi)
    =
    \sum_{p=1}^r
    \left(
      e^{-2\pi i\xi a_p}
      -e^{-2\pi i\xi b_p}
    \right).
\]
Apply Theorem~\ref{thm:indicator-no-aliasing}.  The case \(\xi=0\) follows by continuity.
\end{proof}

\vspace{6pt}

\begin{remark}[Why diameter is the natural quantity]
At frequency \(\xi\), the relevant geometric projection is
\[
    t\longmapsto2\pi\xi t\pmod{2\pi}.
\]
The measure \(T\) controls the total angular length of the projected set,
whereas \(D\) controls whether the projection spans more than one turn.
The hypothesis \(|\xi|D\leq1\) is exactly the condition preventing
multiple wrapping.
\end{remark}

\section{Aliasing and sharpness of the diameter condition}
\label{sec:aliasing}

The preceding theorem cannot be extended by replacing the diameter with
the measure. In fact, the obstruction already occurs for two intervals.

\begin{figure}[ht]
	\centering
	\begin{tikzpicture}[scale=1.05,>=Latex]
		
		\draw[thick] (-0.5,0) -- (8.0,0);
		
		\draw (0,0.08) -- (0,-0.08) node[below=4pt] {$0$};
		\draw (1.6,0.08) -- (1.6,-0.08) node[below=4pt] {$\ell$};
		\draw (5.2,0.08) -- (5.2,-0.04) node[below=4pt] {$\frac1{\xi}$};
		\draw (6.8,0.08) -- (6.8,-0.04) node[below=4pt] {$\frac1{\xi}+\ell$};
		
		\draw[line width=2.2pt,blue] (0,0) -- (1.6,0);
		\draw[line width=2.2pt,blue] (5.2,0) -- (6.8,0);
		
		\node[blue] at (0.8,0.45) {$[0,\ell]$};
		\node[blue] at (6.0,0.45) {$\left[\frac1{\xi},\frac1{\xi}+\ell\right]$};
		
		\draw[<->,thick]
		(0,-0.7) -- (5.1,-0.7)
		node[midway,below=4pt] {$\,\,$ one period $\,=\,\frac1{\xi}$};
		
		\node[above] at (3.4, 0.2) {Real line};
		
		\begin{scope}[shift={(4.2,-4.0)}]
			\def\r{1.8}
			\draw[thick] (0,0) circle (\r);
			
			```
			\draw[thin,gray!60] (-2.2,0) -- (2.2,0);
			\draw[thin,gray!60] (0,-2.2) -- (0,2.2);
			
			\coordinate (A) at ({\r*cos(0)},{\r*sin(0)});
			\coordinate (B) at ({\r*cos(-55)},{\r*sin(-55)});
			
			\draw[line width=3pt,red]
			({\r*cos(0)},{\r*sin(0)})
			arc[start angle=0,end angle=-55,radius=\r];
			
			\node at (0.65,-0.3) {$e^{-2\pi i\xi t}$};
			
			\fill (A) circle (1.4pt);
			\fill (B) circle (1.4pt);
			
			\node[red] at (2.6,-0.95) {same arc};
			
			\node[above right] at (A) {$1$};
			\node[below right] at (B) {$e^{-2\pi i\xi\ell}$};
			
			```
			
		\end{scope}
		
		\node at (1.3,-3.7) {Unit circle};
		
		\draw[->,thick] (0.8,-0.15) .. controls (0.8,-1.5) and (2.0,-4.0) .. (5.4,-5.1);
		\draw[->,thick] (6.0,-0.15) .. controls (6.0,-1.5) and (5.1,-2.4) .. (5.6,-4.7);
		
	\end{tikzpicture}
	\caption{Aliasing in Proposition~4.1. The two intervals
		$[0,\ell]$ and
		$[1/\xi,\,1/\xi+\ell]$
		are separated by one period of the exponential map
		$t\mapsto e^{-2\pi i\xi t}$.
		They therefore project onto the same arc of the unit circle and their
		Fourier contributions at frequency $\xi$ reinforce coherently. In
		particular,
		$|\widehat{\mathbf{1}_\TT}(\xi)|
		=2\sin(\pi\xi\ell)/(\pi\xi)$,
		which is larger than the corresponding Fourier coefficient of the
		centered interval of the same measure.}
	\label{fig:aliasing-prop41}
\end{figure}

\vspace{8pt}


\begin{proposition}[Two-component aliasing]
\label{prop:two-component-aliasing}
Fix \(\xi>0\) and choose \(0<\ell<\frac{1}{2\xi}\). Let
\begin{equation}
    \TT
    =
    [0,\ell]
    \cup
    \left[\frac1\xi,\frac1\xi+\ell\right],
    \label{eq:aliasing-set}
\end{equation}
with size \(|\TT|=2\ell\) and diameter \( 
    \diamess(\TT)=\frac1\xi+\ell
\). 
Then,
\begin{equation}
    |\widehat{\one_\TT}(\xi)|
    =
    \frac{2\sin(\pi\xi\ell)}{\pi\xi}
    >
    \frac{\sin(2\pi\xi\ell)}{\pi\xi}
    =
    \widehat{\one_{[-\ell,\ell]}}(\xi).
    \label{eq:aliasing-reversal}
\end{equation}
\end{proposition}

\begin{proof}
The two components in \eqref{eq:aliasing-set} are separated by exactly
one period of the exponential at frequency \(\xi\).  Therefore their
Fourier contributions coincide:
\[
\begin{aligned}
    \widehat{\one_\TT}(\xi)
    &=
    \left(1+e^{-2\pi i}\right)
    \int_0^\ell e^{-2\pi i\xi t}\,dt\\
    &=
    2e^{-\pi i\xi\ell}
    \frac{\sin(\pi\xi\ell)}{\pi\xi}.
\end{aligned}
\]
The centered interval of the same measure is \([-\ell,\ell]\), and
\[
    \widehat{\one_{[-\ell,\ell]}}(\xi)
    =
    \frac{\sin(2\pi\xi\ell)}{\pi\xi}.
\]
Putting \(x=\pi\xi\ell\in(0,\pi/2)\), it is immediate that
\(
    2\sin x-\sin(2x)
    =
    2\sin x(1-\cos x)>0
\), 
which proves \eqref{eq:aliasing-reversal}.
\end{proof}

\vspace{0pt}

\begin{corollary}[Sharpness of the one-period threshold]
\label{cor:sharp-threshold}
For every \(\varepsilon>0\) and every \(\xi>0\), there exists a
finite union of two intervals \(\TT\) such that
\[
    1<|\xi|\diamess(\TT)<1+\varepsilon
\]
and
\[
    |\widehat{\one_\TT}(\xi)|
    >
    \widehat{\one_{I_T}}(\xi).
\]
\end{corollary}

\begin{proof}
Use Proposition~\ref{prop:two-component-aliasing} and choose
\(
 0<\xi\ell<\min\{1/2,\varepsilon\}
\).
Then \(\xi\diamess(\TT)=1+\xi\ell\).
\end{proof}

\begin{remark}[Measure versus diameter]
In Proposition~\ref{prop:two-component-aliasing}, the product
\(
 \xi|\TT|=2\xi\ell
\)
can be arbitrarily small, while the interval estimate still fails.  Thus
small time--bandwidth product based on measure alone does not prevent
pointwise phase reinforcement between separated components.
\end{remark}
	
\section{Constructive concentration on small-diameter sets}
\label{sec:constructive-concentration}

We now combine the phase-optimization theorem with the pointwise
no-aliasing estimate to obtain our main result.

\begin{lemma}[Quadratic form over a measurable set]
\label{lem:quadratic-form-general}
Let \(\TT\subset\R\) have finite measure and let \(f\in\PW(\WW)\), where \(\WW\) is bounded.  Then
\begin{equation}
    \int_\TT|f(t)|^2\,dt
    =
    \iint_{\WW\times\WW}
       \widehat f(\omega)
       \overline{\widehat f(\eta)}
       \widehat{\one_\TT}(\eta-\omega)
       \,d\omega\,d\eta.
    \label{eq:quadratic-form-general}
\end{equation}
\end{lemma}

\begin{proof}
The proof is the same as that of Lemma~\ref{lem:quadratic-form-interval}, with \(I_T\) replaced by \(\TT\). The boundedness of \(\WW\) gives
\(
 \widehat f\in L^1(\WW)
\),
and the finite measure of \(\TT\) justifies the application of Fubini's theorem.
\end{proof}

We can now prove Theorem~\ref{thm:main-introduction}. The preceding lemma puts the concentration on \(\TT\) in a form where the two ingredients developed in Sections~\ref{sec:phase-optimization} and \ref{sec:no-aliasing} can be combined directly. Indeed, if \(\WW\) is an interval of length \(W\), then every difference frequency \(\xi=\eta-\omega\), with \(\omega,\eta\in\WW\), satisfies \(|\xi|\leq W\). Hence the hypothesis
\[
W\diamess(\TT)\leq1
\]
places every such frequency in the no-aliasing regime, so that the Fourier coefficient of \(\one_\TT\) is dominated in modulus by that of the centered interval \(I_T  \). At the same time, since  \(T=|\TT|\leq\diamess(\TT)\), the same hypothesis implies \(WT\leq1\), and the Fourier transform of \(\one_{I_T}\) is therefore nonnegative throughout the whole difference set \(\WW-\WW\). These two facts allow us first to replace the kernel associated with \(\TT\) by the centered-interval kernel and then to align the phases of \(\widehat f\). The resulting quadratic form is precisely the concentration on \(I_T\) of the function \(g\) defined by \(\widehat g=|\widehat f|\).

\begin{proof}[Proof of Theorem~\ref{thm:main-introduction}]
Since \(\WW\) is an interval of length \(W\), it holds \(|\eta-\omega|\leq W\) for \(\omega,\eta\in\WW\). Therefore
\[
    |\eta-\omega|D\leq WD\leq1.
\]
Theorem~\ref{thm:indicator-no-aliasing} gives
\begin{equation}
    |\widehat{\one_\TT}(\eta-\omega)|
    \leq
    \widehat{\one_{I_T}}(\eta-\omega).
    \label{eq:kernel-domination-main-proof}
\end{equation}
Since \(T\leq D\), condition \eqref{eq:main-condition-intro} implies \(WT\leq1\), and the interval kernel in \eqref{eq:kernel-domination-main-proof} is nonnegative on the difference set \(\WW - \WW\).

Using Lemma~\ref{lem:quadratic-form-general},
\[
\begin{aligned}
    \int_\TT|f(t)|^2\,dt
    &\leq
    \iint_{\WW\times\WW}
       |\hat f(\omega)|\,|\hat f(\eta)|
       |\widehat{\one_\TT}(\eta-\omega)|
       \,d\omega\,d\eta\\
    &\leq
    \iint_{\WW\times\WW}
       |\hat f(\omega)|\,|\hat f(\eta)|
       \widehat{\one_{I_T}}(\eta-\omega)
       \,d\omega\,d\eta\\
    &=
    \int_{I_T}|g(t)|^2\,dt,
\end{aligned}
\]
where the last identity follows from Lemma~\ref{lem:quadratic-form-interval} and the definition of $g$ in \eqref{eq:g-main-section}.  Plancherel's theorem gives
\(
 \|g\|_2=\|f\|_2
\),
which proves \eqref{eq:normalized-function-by-function}.
\end{proof}

\begin{corollary}[A no-aliasing Donoho--Stark inequality]
\label{cor:local-operator}
Under the hypotheses of Theorem~\ref{thm:main-introduction},
\begin{equation}
    \cC_\WW(\TT)
    \leq
    \cC_\WW(I_T).
    \label{eq:local-operator}
\end{equation}
\end{corollary}

\begin{proof}
For each nonzero \(f\in\PW(\WW)\), Theorem \ref{thm:main-introduction} supplies a function \(g\in\PW(\WW)\) with the same norm such that
\[
    c_\TT(f)
    \leq
    c_{I_T}(g)
    \leq
    \cC_\WW(I_T).
\]
Taking the supremum over \(f\) yields the result.
\end{proof}

\begin{remark}[Invariances]
The statement is invariant under translations of \(\TT\) and translations
of the frequency interval \(\WW\). indeed, translating \(\TT\) multiplies
\(\widehat{\one_\TT}\) by a unimodular phase, while translating \(\WW\)
corresponds to modulation in the time variable and leaves \(|f|\)
unchanged.
\end{remark}

\section{Discrete version on the unit circle}
\label{sec:polynomial}

We now introduce a discrete analogue of the previous concentration problem.
As in the original version of this argument, the continuous Fourier support is
replaced by the degree of a polynomial, while the observation set is now a
measurable subset of the unit circle.

Let \(
    \Torus=\{z\in\mathbb C:|z|=1\}
    \simeq \R/(2\pi\mathbb Z)
\) be the unit circle. We denote by \(\dm(z)\) the Lebesgue (arc-length) measure on \(\Torus\),
normalized so that \(
    \int_{\Torus}\dm(z)=2\pi
\). More explicitly, if \(\Omega\subseteq\Torus\) is measurable and
\[
    \Omega=\{e^{i\theta}:\theta\in\Theta\},
    \qquad
    \Theta\subset[-\pi,\pi)
\]
up to a null set, then
\begin{equation}
    |\Omega|_{\Torus}
    :=
    \int_\Omega\dm(z)
    =
    \int_\Theta d\theta,
    \label{eq:circle-measure}
\end{equation}
and through this section we write \( T=|\Omega|_{\Torus}\). The centered arc of the same measure \(T\) will be denoted by
\begin{equation}
    I_T
    =
    \{e^{i\theta}: -T/2 \leq \theta \leq T/2\}.
    \label{eq:centered-circle-arc}
\end{equation}
For a polynomial
\(
    P(z)=a_0+a_1z+\cdots+a_nz^n
    \in\mathbb C_n[z] = \{R\in \C[z] \colon \deg R \leq n \},
\)
we define its concentration on \(\Omega\) by
\begin{equation}
    c_\Omega(P)
    =
    \frac{\displaystyle\int_\Omega |P(z)|^2\,\dm(z)}
         {\displaystyle\int_\Torus |P(z)|^2\,\dm(z)}.
    \label{eq:circle-concentration}
\end{equation}
By Parseval's identity,
\begin{equation}
    \|P\|_{L^2(\Torus)}^2
    =
    \int_\Torus |P(z)|^2\,\dm(z)
    =
    \int_{-\pi}^{\pi}|P(e^{i\theta})|^2\,d\theta
    =
    2\pi\sum_{k=0}^n|a_k|^2.
    \label{eq:circle-polynomial-norm}
\end{equation}
It is also convenient to introduce the optimal polynomial concentration
\[
    \mathcal C_n(\Omega)
    =
    \sup_{0\neq P\in\mathbb C_n[z]}c_\Omega(P).
\]
The discrete analogue of the Donoho--Stark concentration conjecture can then
be formulated as follows.

\begin{conjecture}[Discrete concentration conjecture]
\label{conj:discrete}
Fix \(n\geq1\) and \(0<T\leq2\pi\). Among all measurable sets
\(\Omega\subseteq\Torus\) satisfying \(|\Omega|_{\Torus}=T\), the quantity \(\mathcal C_n(\Omega)\) is maximized by an arc of length
\(T\). Equivalently,
\[
    \mathcal C_n(\Omega)
    \leq
    \mathcal C_n(I_T).
\]
\end{conjecture}

As in the continuous case, the position of the comparison interval is
irrelevant. Indeed, if
\[
    J=e^{i\theta_0}I_T
\]
for some $\theta_0$, and
\[
    P(z)=\sum_{k=0}^n a_kz^k,
    \qquad
    P_{\theta_0}(z)=\sum_{k=0}^n a_ke^{-ik\theta_0}z^k,
\]
then
\[
    c_{I_T}(P)=c_J(P_{\theta_0}).
\]
We may therefore work throughout with the centered arc \(I_T\).

The correction introduced in the continuous problem also has a direct
counterpart on the circle. The measure \(T\) alone is not sufficient to
prevent different components of \(\Omega\) from becoming identified by a
higher Fourier mode. We therefore need the circular analogue of the
essential diameter.

\begin{definition}[Essential arc diameter]
\label{def:essential-arc-diameter}
For a measurable set \(\Omega\subseteq\Torus\), define
\begin{equation}
    D_{\Torus}(\Omega)
    =
    \inf\left\{
       D\in[0,2\pi]:
       \Omega\text{ is contained a.e. in a circular arc of length }D
    \right\}.
    \label{eq:essential-arc-diameter}
\end{equation}
\end{definition}

The following lemma is the discrete counterpart of Theorem~\ref{thm:indicator-no-aliasing}, where here the integer \(q\) plays the role of a difference frequency.

\vspace{6pt}

\begin{lemma}[Fourier coefficients in the no-aliasing regime]
\label{lem:circle-fourier-coefficient}
Let \(\Omega\subseteq\Torus\) be measurable, with positive finite measure \(    T=|\Omega|_{\Torus}\), and finite diameter \(D=D_{\Torus}(\Omega)\), and let \(q\geq1\) be an integer. If
\begin{equation}
    qD\leq2\pi,
    \label{eq:qD-circle}
\end{equation}
then
\begin{equation}
    \left|
       \int_\Omega z^q\,\dm(z)
    \right|
    =
    \left|
       \int_\Theta e^{iq\theta}\,d\theta
    \right|
    \leq
    \frac{2\sin(qT/2)}{q}.
    \label{eq:circle-fourier-coefficient}
\end{equation}
For \(q=0\), the corresponding value is \(T\).
\end{lemma}

\begin{proof}
Up to a rotation and a null set, choose a lift
\[
    \widetilde\Theta\subset[a,a+D]\subset\R.
\]
Under the change of variables \(s=q\theta\), the set
\(q\widetilde\Theta\) is contained in an interval of length
\(qD\leq2\pi\) and has measure \(qT\). Hence
Lemma~\ref{lem:single-arc} gives
\[
\begin{aligned}
    \left|
       \int_\Theta e^{iq\theta}\,d\theta
    \right|
    &=
    \frac1q
    \left|
       \int_{q\widetilde\Theta}e^{is}\,ds
    \right|\\
    &\leq
    \frac{2\sin(qT/2)}q.
\end{aligned}
\]
Since \(T\leq D\), condition \eqref{eq:qD-circle} also implies \(    0\leq  {qT}/{2}\leq\pi\), so the right hand side is indeed non-negative.
\end{proof}

We are now ready to state the discrete version of the constructive
concentration theorem. Given
\[
    P(z)=\sum_{k=0}^n a_kz^k,
\]
define
\begin{equation}
    Q(z)=\sum_{k=0}^n|a_k|z^k.
    \label{eq:polynomial-phase-aligned}
\end{equation}
Thus \(Q\) is obtained from \(P\) by removing the phases of its
coefficients, exactly as \(\widehat g=|\widehat f|\) removes the Fourier
phase in the continuous setting.

\vspace{10pt}

\begin{theorem}[Constructive polynomial concentration]
\label{thm:polynomial-phase-symmetrization}
Let \(\Omega\subseteq\Torus\) be measurable, with positive finite measure \(    T=|\Omega|_{\Torus}\), and finite diameter \(D=D_{\Torus}(\Omega)\), and let \(
   P \in \C_n[z]
\) be a polynomial of degree \(n\geq1\). Assume
\begin{equation}
    nD\leq2\pi.
    \label{eq:polynomial-diameter-condition}
\end{equation}
Then, for \(Q\) defined in \eqref{eq:polynomial-phase-aligned}, it holds
\begin{equation}
    \int_\Omega |P(z)|^2\,\dm(z)
    \leq
    \int_{I_T}|Q(z)|^2\,\dm(z).
    \label{eq:polynomial-phase-symmetrization}
\end{equation}
Moreover, \(\|P\|_{L^2(\Torus)} = \|Q\|_{L^2(\Torus)}\), and hence
\begin{equation}
    c_\Omega(P)\leq c_{I_T}(Q).
    \label{eq:polynomial-concentration-normalized}
\end{equation}
\end{theorem}

\begin{proof}
Expanding the square in angular coordinates,
\begin{equation}
\begin{aligned}
    \int_\Omega|P(z)|^2\,\dm(z)
    &=
    \sum_{l,m=0}^n
       a_l\overline{a_m}
       \int_\Theta e^{i(l-m)\theta}\,d\theta.
\end{aligned}
\label{eq:polynomial-quadratic-expansion}
\end{equation}
For \(q=|l-m|\), one has \(q\leq n\), and therefore
\[
    qD\leq nD\leq2\pi.
\]
By Lemma~\ref{lem:circle-fourier-coefficient},
\[
    \left|
       \int_\Theta e^{i(l-m)\theta}\,d\theta
    \right|
    \leq s_{|l-m|},
\]
where
\begin{equation}
    s_0=T,
    \qquad
    s_q=\frac{2\sin(qT/2)}q,
    \quad 1\leq q\leq n.
    \label{eq:centered-arc-toeplitz}
\end{equation}
Since \(T\leq D\), hypothesis
\eqref{eq:polynomial-diameter-condition} implies \( nT/2 \leq\pi \), and hence all the coefficients \(s_q\) are nonnegative. Thus
\[
\begin{aligned}
    \int_\Omega|P(z)|^2\,\dm(z)
    &\leq
    \sum_{l,m=0}^n
       |a_l|\,|a_m|\,s_{|l-m|}\\
    &=
    \int_{I_T}|Q(z)|^2\,\dm(z).
\end{aligned}
\]
Finally, by \eqref{eq:circle-polynomial-norm}, we have \(\|P\|_{L^2(\Torus)} =     \|Q\|_{L^2(\Torus)} \), which gives \eqref{eq:polynomial-concentration-normalized}.
\end{proof}

As an immediate consequence, we also have a positive result for Conjecture~\ref{conj:discrete} whenever the underlying set satisfies the no-aliasing condition \(nD_{\Torus}(\Omega)\leq2\pi\).

\begin{corollary}[Discrete concentration in the no-aliasing regime]
	\label{cor:discrete-concentration}
	Let \(n\geq1\), and let \(\Omega\subseteq\Torus\) be measurable with
	measure \(T=|\Omega|_{\Torus}\) and diameter \(D=D_{\Torus}(\Omega)\). If
	\begin{equation}
		nD\leq2\pi,
		\label{eq:discrete-concentration-diameter}
	\end{equation}
	then
	\begin{equation}
		\mathcal C_n(\Omega)
		\leq
		\mathcal C_n(I_T).
		\label{eq:discrete-concentration-corollary}
	\end{equation}
\end{corollary}

\begin{proof}
	Let \(0\neq P\in\mathbb C_n[z]\), and let \(m\leq n\) be its degree.
	Since
	\[
	mD\leq nD\leq2\pi,
	\]
	Theorem~\ref{thm:polynomial-phase-symmetrization} gives a polynomial
	\(Q\in\mathbb C_n[z]\) such that
	\[
	c_\Omega(P)\leq c_{I_T}(Q)
	\leq \mathcal C_n(I_T).
	\]
	Taking the supremum over all nonzero \(P\in\mathbb C_n[z]\) yields the result.
\end{proof}

Notice that the old condition involving only the measure,
\[
    \frac{nT}{2}\leq\pi,
\]
is sufficient to make the centered-arc kernel positive, but it is not
sufficient to justify the comparison with an arbitrary disconnected
\(\Omega\). The stronger condition
\[
    nD_{\Torus}(\Omega)\leq2\pi
\]
is precisely the no-aliasing condition needed to control every difference
frequency \(0\leq q\leq n\), therefore completely parallel to the continuous condition
\(
    W\,\diamess(\TT)\leq1
\).
The notation is intentionally only partly shared: \(T\) and \(I_T\) denote
the size of the observation set and the centered comparator in both
settings, while \(\TT\subset\R\) and \(\Omega\subset\Torus\) distinguish
the two ambient spaces.

\section{Coefficient rearrangement and Montgomery's inequality}
\label{sec:coefficient-rearrangement}

The polynomial obtained in Theorem~\ref{thm:polynomial-phase-symmetrization}
has positive coefficients, but their positions have not yet been
rearranged. We now study this second discrete rearrangement and relate it
explicitly to a classical result of Montgomery.

\subsection{Montgomery's result in the present notation}

In \cite{montgomery1976note}, Montgomery proved a rearrangement inequality
for a more general family of uniformly bounded functions satisfying a
Bessel-type inequality. Restricting his result to the circle
\(\Torus=\R/(2\pi\mathbb Z)\) and to the functions
\(\varphi_k(\theta)=\cos(k\theta)\), it takes the following form.

\begin{theorem}[Montgomery, Thm.~1~\cite{montgomery1976note}]
Let \(
    f(\theta)=\sum_{k=0}^{\infty}b_k\cos(k\theta),
\)
and define
\[
    f^{**}(\theta)
    =
    \sum_{k=0}^{\infty}b_k^{**}\cos(k\theta),
\]
where the numbers \(b_k^{**}\) are the values \(|b_k|\), rearranged so
that \(
    b_0^{**}\geq b_1^{**}\geq b_2^{**}\geq\cdots
\).
Then, for every measurable \(\Omega\subseteq\Torus\), with measure \(
    T=|\Omega|_{\Torus},
\)
one has
\begin{equation}
    \int_\Omega |f(\theta)|^2\,d\theta
    \leq
    20\int_{I_T}|f^{**}(\theta)|^2\,d\theta.
    \label{eq:montgomery-cosine}
\end{equation}
\end{theorem}

The rearrangement \(f\mapsto f^{**}\) may be regarded as a discrete
counterpart of the symmetric decreasing rearrangement appearing in the
Donoho--Stark argument: the magnitudes of the coefficients are placed in
decreasing order as the frequency moves away from the origin.

For a finite cosine polynomial
\[
    f(\theta)=\sum_{k=0}^{N}b_k\cos(k\theta),
\]
we may rewrite
\begin{equation}
\begin{aligned}
    2e^{iN\theta}f(\theta)
    &=
    b_N+b_{N-1}e^{i\theta}+\cdots+b_1e^{i(N-1)\theta}
    +2b_0e^{iN\theta}\\
    &\qquad
    +b_1e^{i(N+1)\theta}+\cdots
    +b_{N-1}e^{i(2N-1)\theta}
    +b_Ne^{i2N\theta}.
\end{aligned}
\label{eq:cosine-to-polynomial}
\end{equation}
Thus, setting \(z=e^{i\theta}\), it is natural to associate with \(f\)
the symmetric polynomial
\begin{equation}
    P(z)
    =
    b_N+b_{N-1}z+\cdots+b_1z^{N-1} +b_0z^N +b_1z^{N+1}+\cdots+b_Nz^{2N},
    \label{eq:montgomery-symmetric-polynomial}
\end{equation}
for which \(|P(e^{i\theta})|=2|f(\theta)|\).

If \(b_0^{**},\ldots,b_N^{**}\) are the numbers \(|b_0|,\ldots,|b_N|\) arranged so that \(
    b_0^{**}\geq b_1^{**}\geq\cdots\geq b_N^{**},
\)
and
\begin{equation}
    P^{**}(z)
    =
    b_N^{**}+b_{N-1}^{**}z+\cdots+b_1^{**}z^{N-1}
    +2b_0^{**}z^N
    +b_1^{**}z^{N+1}+\cdots+b_N^{**}z^{2N},
    \label{eq:montgomery-polynomial-rearrangement}
\end{equation}
then Montgomery's theorem becomes, in our polynomial notation,
\begin{equation}
    \int_\Omega|P(z)|^2\,\dm(z)
    \leq
    20\int_{I_T}|P^{**}(z)|^2\,\dm(z).
    \label{eq:montgomery-polynomial}
\end{equation}
Hence Montgomery obtains a universal factor \(20\), with no diameter restriction, for this symmetric class of even-degree polynomials.

For comparison, a constant-one rearrangement over the whole unit circle
appears in the classical result of Gabriel.

\begin{theorem}[Gabriel, Thm.~4~\cite{gabriel1932rearrangement}]
Let \(k\geq1\) be an integer and
\[
    P(\theta)=\sum_{r=-R}^{R}a_re^{ir\theta},
    \qquad
    P^+(\theta)=\sum_{r=-R}^{R}a_r^+e^{ir\theta},
\]
where the \(a_r^+\) are the numbers \(|a_r|\) ordered so that \(
    a_0^+\geq a_{-1}^+\geq a_1^+
    \geq a_{-2}^+\geq a_2^+\geq\cdots
\).
Then
\[
    \int_0^{2\pi}|P(\theta)|^{2k}\,d\theta
    \leq
    \int_0^{2\pi}|P^+(\theta)|^{2k}\,d\theta.
\]
\end{theorem}

It should be mentioned that this is actually an improvement on the same result by Hardy and Littlewood appearing in \cite{hardy1948new}, but which had some additional symmetry hypothesis on the coefficients.

\subsection{Rearrangement of the polynomial coefficients}

We now return to the centered arc \(I_T\). After the phase removal of Section~\ref{sec:polynomial}, it remains to determine how the positive coefficient magnitudes should be positioned. The relevant result is the following bilinear rearrangement theorem.

\begin{lemma}[Hardy, Littlewood and P\'olya, Thm.~1~\cite{hardy1926maximum}]
\label{lem:discrete-toeplitz-rearrangement}
Consider the bilinear form
\begin{equation}
    S(x,y)
    =
    \sum_{l=0}^{n}\sum_{m=0}^{n}
       s_{l-m}x_ly_m.
    \label{eq:HLP-bilinear-form}
\end{equation}
Suppose that the given coefficients satisfy \(s_\nu=s_{-\nu}\) and \(
    s_0\geq s_1\geq\cdots\geq s_n\geq0
\),
and that the variables \(x_l \geq0\), \(y_m \geq 0\) are given in every respect except their arrangement. Then, among the arrangements for which \(S(x,y)\) assumes its maximum value, there is one for which:

\begin{enumerate}[label=\textup{(\roman*)}]
\item
\(
    x_\mu\leq x_{\mu'}
    \text{  if  }
    \left|\mu'- n/ 2\right|
    <
    \left|\mu- n/ 2\right|
\)

\item no two of \(x_\mu - x_{\mu '}\), where \( \mu < \mu' \), \(|\mu' - n/2| = |\mu - n/2| \), have different signs.
\end{enumerate}

and analogous conditions for variables \(y_\eta\),
\end{lemma}

In the present quadratic case we take \(x_l=y_l\), so one maximizing arrangement places the largest coefficient at the center and then alternates consistently to the right and to the left. More explicitly, if
\[
    P(z)=\sum_{k=0}^{n}a_kz^k,
\]
we denote by \(
    a_0^*,a_1^*,\ldots,a_n^*
\)
a permutation of
\(
    |a_0|,|a_1|,\ldots,|a_n|
\)
which may be chosen to satisfy

\begin{equation}
\begin{cases}
    a_{\lceil n/2\rceil}^*
    \geq
    a_{\lfloor n/2\rfloor}^*
    \geq
    a_{\lceil n/2\rceil+1}^*
    \geq
    a_{\lfloor n/2\rfloor-1}^*
    \geq\cdots
    & n\ \text{odd},\\[2mm]
    a_{n/2}^*
    \geq
    a_{n/2+1}^*
    \geq
    a_{n/2-1}^*
    \geq
    a_{n/2+2}^*
    \geq
    a_{n/2-2}^*
    \geq\cdots
    & n\ \text{even}.
\end{cases}
\label{eq:explicit-coefficient-order}
\end{equation}
Define
\begin{equation}
    P^*(z)=\sum_{k=0}^{n}a_k^*z^k.
    \label{eq:P-star}
\end{equation}
The centered-arc quadratic form is precisely of the type appearing in Lemma~\ref{lem:discrete-toeplitz-rearrangement}. Indeed, with
\begin{equation}
    s_0=T,
    \qquad
    s_\nu=
    \frac{2\sin(\nu T/2)}{\nu},
    \quad 1\leq \nu\leq n,
    \qquad
    s_{-\nu}=s_\nu,
    \label{eq:rearrangement-toeplitz-weights}
\end{equation}
one has
\begin{equation}
    \int_{I_T}
       \left|\sum_{k=0}^{n}a_kz^k\right|^2
       \dm(z)
    =
    \sum_{l,m=0}^{n}
       s_{l-m}a_l\overline{a_m}.
    \label{eq:centered-arc-quadratic-form}
\end{equation}
If \(nT/2\leq\pi\), then \(s_0\geq s_1\geq\cdots\geq s_n\geq0\), because 
\(
    s_\nu
    =
    T\,
    \frac{\sin(\nu T/2)}{\nu T/2}
\)
and \(x\mapsto\sin x/x\) is nonnegative and decreasing on
\([0,\pi]\).

We can now combine the phase removal of Theorem~\ref{thm:polynomial-phase-symmetrization} with this discrete coefficient rearrangement.

\begin{theorem}[Phase and coefficient rearrangement]
\label{thm:phase-coefficient-rearrangement}
Let \(\Omega\subseteq\Torus\) be measurable, with positive finite measure \(    T=|\Omega|_{\Torus}\), and finite diameter \(D=D_{\Torus}(\Omega)\), and let \(
P \in \C_n[z]
\) be a polynomial of degree \(n\geq1\). Assume
\begin{equation}
	nD\leq2\pi.
    \label{eq:coefficient-rearrangement-diameter}
\end{equation}
Let \(P^*\) be obtained by first replacing the coefficients of \(P\) by their absolute values and then arranging them as in \eqref{eq:explicit-coefficient-order}. Then
\begin{equation}
    \int_\Omega|P(z)|^2\,\dm(z)
    \leq
    \int_{I_T}|P^*(z)|^2\,\dm(z).
    \label{eq:phase-coefficient-rearrangement}
\end{equation}
In particular,
\begin{equation}
    c_\Omega(P)\leq c_{I_T}(P^*).
    \label{eq:phase-coefficient-concentration}
\end{equation}
\end{theorem}

\begin{proof}
By Theorem~\ref{thm:polynomial-phase-symmetrization}, using the hypothesis \(nD \leq 2 \pi\), we have
\[
    \int_\Omega|P(z)|^2\,\dm(z)
    \leq
    \int_{I_T}
       \left|
          \sum_{k=0}^{n}|a_k|z^k
       \right|^2
       \dm(z).
\]
Since \(T\leq D\), hypothesis
\eqref{eq:coefficient-rearrangement-diameter} implies \({nT}/{2}\leq\pi\), and therefore the coefficients \(s_\nu\) in \eqref{eq:rearrangement-toeplitz-weights} satisfy all the hypotheses of Lemma~\ref{lem:discrete-toeplitz-rearrangement}. Taking \(x_l=y_l=|a_l| \) and applying that lemma gives
\[
\begin{aligned}
    \int_{I_T}
       \left|
          \sum_{k=0}^{n}|a_k|z^k
       \right|^2
       \dm(z)
    &=
    \sum_{l,m=0}^{n}s_{l-m}|a_l||a_m|\\
    &\leq
    \sum_{l,m=0}^{n}s_{l-m}a_l^*a_m^*\\
    &=
    \int_{I_T}|P^*(z)|^2\,\dm(z).
\end{aligned}
\]
Combining the two inequalities proves \eqref{eq:phase-coefficient-rearrangement}. Finally, \(P\) and \(P^*\) have the same multiset of coefficient magnitudes, so Parseval's identity gives \(\|P\|_{L^2(\Torus)} = \|P^*\|_{L^2(\Torus)}\), and hence \eqref{eq:phase-coefficient-concentration}.
\end{proof}

The relation with Montgomery's theorem can now be stated precisely. For a symmetric polynomial of the form \eqref{eq:montgomery-symmetric-polynomial}, of degree at most \(2N\), Theorem~\ref{thm:phase-coefficient-rearrangement} applies whenever
\[
2N D_{\Torus}(\Omega)\leq 2\pi.
\]
In this no-aliasing regime it gives a constant-one comparison with the rearranged polynomial \(P^*\). Notice, however, that \(P^*\) is in general different from Montgomery's \(P^{**}\). Indeed, the former rearranges the coefficient magnitudes of the associated polynomial (including the extra $2$ factor in the $0$-th coefficient), whereas the latter is induced by rearranging the cosine coefficients before passing to the symmetric polynomial. Thus the present result should be viewed as an analogous constant-one rearrangement inequality with a different explicit comparator, valid under the additional no-aliasing condition on the diameter.

\section{Scope and limitations of the method}
\label{sec:limitations}

The results above give explicit function-by-function comparisons in the no-aliasing regime. They do not, however, resolve the unrestricted Donoho--Stark concentration conjecture. We close by explaining the difference between the pointwise comparison used in the proof and the operator-level statement of the conjecture.

The distinction is particularly clear at the operator level. Recall that \(P_\TT\) and \(P_\WW\) denote, respectively, the time- and frequency-limiting orthogonal projections. The positive concentration operator on \(\PW(\WW)\) is
\[
    K_{\TT,\WW}
    =
    \left.P_\WW P_\TT P_\WW\right|_{\PW(\WW)},
    \qquad
    \mathcal C_\WW(\TT)=\|K_{\TT,\WW}\|.
\]
Thus the unrestricted conjecture compares the norms of the full concentration operators, whereas the present proof proceeds through a stronger frequency-by-frequency estimate. Indeed, it establishes for \(T=|\TT|\) and \(I_T=[-T/2,T/2]\) that
\begin{equation}
    |\widehat{\one_\TT}(\xi)|
    \leq
    \widehat{\one_{I_T}}(\xi)
    \label{eq:stronger-modewise}
\end{equation}
for every difference frequency \(\xi \in \WW - \WW\) occurring in the quadratic form, provided that
\[
    |\xi|\diamess(\TT)\leq1.
\]
Once \eqref{eq:stronger-modewise} is available, the function-by-function concentration comparison follows by replacing the spectral factors by their absolute values. Proposition \ref{prop:two-component-aliasing} and Corollary \ref{cor:sharp-threshold}, however, show that this pointwise estimate fails for general disconnected sets as soon as the one-period condition is lost, even when the product \(|\xi|T\) is arbitrarily small.

This failure does not by itself contradict the operator-level conjecture. Indeed, for \(f\in\PW(\WW)\),
\begin{equation}
    \int_\TT|f(t)|^2\,dt
    =
    \iint_{\WW\times\WW}
       \widehat f(\omega)
       \overline{\widehat f(\eta)}
       \widehat{\one_\TT}(\eta-\omega)
       \,d\omega\,d\eta,
    \label{eq:full-quadratic-form-limitations}
\end{equation}
and the value of the quadratic form therefore depends on all difference frequencies simultaneously, coupled through the spectral factors of \(f\). A pointwise failure of modulus domination at one frequency need not determine the ordering of the corresponding operator norms.

The role of the disconnected components can be seen more geometrically. Suppose, for simplicity, that
\[
    \TT=\bigcup_j J_j,
    \qquad
    J_j=x_j+I_j,
\]
where the \(I_j\) are centered intervals, hence
\begin{equation}
    \widehat{\one_\TT}(\xi)
    =
    \sum_j
       e^{-2\pi i x_j\xi}
       \widehat{\one_{I_j}}(\xi).
    \label{eq:component-interaction-final}
\end{equation}
For a fixed frequency \(\xi\), the factors
\[
    e^{-2\pi i x_j\xi}
\]
may be viewed as points on the unit circle. Thus the relative positions of the components on the real line are seen only modulo the period \(1/|\xi|\). Components that are far apart may therefore have identical phases and reinforce coherently, as in the two-component example of Section~\ref{sec:aliasing}.

The important point is that this phase configuration depends on \(\xi\). The same collection of component centers \(\{x_j\}\) produces a different configuration on the unit circle for each difference frequency. Hence there is no single phase arrangement of the components that controls the whole concentration problem. The full quadratic form \eqref{eq:full-quadratic-form-limitations} combines these frequency-dependent interactions simultaneously. Equivalently, the full concentration operator retains the interactions between different components of \(\TT\), rather than controlling their contributions independently at each Fourier mode.

The limitation encountered here has a direct parallel in the rearrangement argument of Donoho and Stark. Both approaches give explicit function-by-function comparisons only in restricted regimes: symmetric decreasing rearrangement in their time-limited formulation, and Fourier-phase alignment together with the no-aliasing estimate in the present band-limited formulation. Since these prescribed comparisons are stronger than the corresponding operator-level extremal statements, their failure outside the respective regimes does not imply failure of the conjecture itself.

Taken together, the present results identify a precise regime in which a constructive concentration comparison is available: phase alignment gives the optimal centered-interval comparator for a prescribed Fourier modulus, while the one-period condition prevents geometric aliasing of the observation set. The same mechanism yields corresponding polynomial and coefficient-rearrangement results on the circle. The two-component examples show, however, that aliasing obstructs this frequency-by-frequency method beyond the one-period regime. This is an obstruction to the present method, not evidence against the unrestricted concentration conjecture; any extension beyond this regime must therefore exploit interactions between different frequencies, or some other genuinely global mechanism.

\newpage
\bibliographystyle{unsrt}
\bibliography{references}

\end{document}